\documentclass{article}

\usepackage[english]{babel}

\usepackage[letterpaper,top=2cm,bottom=2cm,left=3cm,right=3cm,marginparwidth=1.75cm]{geometry}

\usepackage{amsmath}
\usepackage{graphicx}
\usepackage{dirtytalk}
\usepackage{authblk}
\usepackage[colorlinks=true, allcolors=blue]{hyperref}

\title{An attempt to trace the birth of importance sampling}
\author{Charly Andral\thanks{The author is supported by a grant from Région Ile-de-France}\\ andral@ceremade.dauphine.fr}

\affil{CEREMADE, CNRS, UMR 7534, Université Paris-Dauphine, PSL University, 75016 Paris, France }

\begin{document}
\maketitle

\begin{abstract}
    In this note, we try to trace the birth of importance sampling (IS) back to 1949. We found the classical formulation of IS in a paper from Kahn in June 1949. As for the appearance of the expression \textit{importance sampling} itself, it may have appeared a few months later, maybe in 1949 during conferences, but there is no published article with the name before 1950.
\end{abstract}
\section{Context}

Importance sampling is now well known but when compared to other methods (e.g. Metropolis algorithm), its birthdate is difficult to establish. Two researchers in the United States seem to have had a role in the birth of importance sampling : Herman Kahn, Gerald Goertzel. They are part of a larger group of scientists who worked on the problem of neutron and gamma
ray attenuation in thick shields. It includes J. von Neumann, S. Ulam., H. Bethe, W. DeMarcus, R. Echols, U. Fano, T.E. Harris, A.S, Householder
G. Goertzel, N. Metropolis, L.A. Ohlinger, R.D, Richtmyer,
and S.W.W. Shor \cite{kahnStochasticMonteCarlo1949}. 

\section{Chronology}

\subsection{1949: Quota Sampling}

The true chronology is difficult to establish as papers exist in different versions and with different dates and as the thinking and writing of the corresponding articles may have occurred much earlier than the publication date. No article mentions work on importance sampling before 1949. 

The first mention of the idea of importance sampling is found in two papers of June 1949, one from Kahn \cite{kahnStochasticMonteCarlo1949a}(June 14), and the other from  Goertzel \cite{goertzelQuotaSamplingImportance1949} (June 21) under the name \textit{quota sampling}. Note that \cite{kahnStochasticMonteCarlo1949a} exists under another reference from the RAND Corporation (where Kahn worked) with a revised version \cite{kahnStochasticMonteCarlo1949} with the date July 14, 1949. Goertzel's article is more focused on the physical applications of the method while \cite{kahnStochasticMonteCarlo1949a} contains already the classical justification of importance sampling
 \begin{equation*}
    J = \int \frac{g(x)f(x)}{f^*(x)}f^*(x)dx,
\end{equation*}
where  $f$ is the target density, $g$ the function to integrate, and $f^*$ an arbitrary density function, meaning that one can sample according to $f^*$ to estimate $J$. 

However, no article published in 1949 mentions \textit{importance sampling}. The term \textit{importance} only appears in \cite{goertzelQuotaSamplingImportance1949} as \textit{importance function}.

\subsection{1950: Importance Sampling}

At least three 1950 papers use the term \textit{importance sampling} as what was previously known under the name \textit{quota sampling}. The first two  are in the proceedings of an IBM conference hold in New York in November 1949, published in 1950 \cite{internationalbusinessmachinescorporationProceedingsSeminarScientific1950}.
One is by Kahn (presenting work done with T.E. Harris) \cite{kahnModificationMonteCarlo1950} and has a very informative paragraph about \textit{importance sampling} :

\begin{quote}
\begin{minipage}{\linewidth}
    Generalizing a little, we can say that in making
    the $f(x) \to f^*(x)$ transformation we are trying to
    sample every portion of the phase space that we are
    studying, with a frequency proportional to the importance of that part of phase space times the probability of getting into it. For this reason the name
    \say{importance sampling}\footnote[1]{We had previously called this process \say{quota sampling.} This last terminology is a little misleading as it is not identical with the quota sampling of the statisticians.} has been suggested by
    G. Goertzel for this procedure. In this terminology $g(x)$ is called the importance function.
    In general the importance of any region is directly proportional
    to the amount that the region contributes to the answer.
\end{minipage}
\end{quote}
The second article in \cite{internationalbusinessmachinescorporationProceedingsSeminarScientific1950} is from J.H. Curtiss \cite{curtissSamplingMethodsApplied1950} and contains the following footnote about \textit{importance sampling}
\begin{quote}
See also the paper by Kahn and Harris. The
method is called that of \say{importance sampling}
by these writers. Dr. W. Edwards Deming pointed out to the author that the method has been known
for many years in various forms to sampling experts working in the social sciences, who have sometimes called it \say{sampling with probability
in proportion to size}.
\end{quote}
At this point, the chronology starts to get fuzzy. The Kahn and Harris paper has the following reference : 
\begin{quote}
    The Monte Carlo Method: Proceedings of a Symposium held in Los Angeles, Calif., on June 29, 30 and July 1, 1949. To be published by the National Bureau of Standards.
\end{quote}
...but it was in fact published in 1951 \cite{kahnEstimationParticleTransmission1951} in \cite{analysisu.s.MonteCarloMethod1951}. 

The third 1950 paper \cite{goertzelMonteCarloMethods1950} is dated from February 1950, where Goertzel and Kahn write this footnote:

\begin{quote}
    
Importance sampling is the term used herein for what has in several reports been called \say{quota sampling.} \say{Importance sampling} is considered to be a more descriptive term.
\end{quote}

\subsection{1951: a 1949 paper?}
The last paper we will mention here is the Kahn and Karris paper \cite{kahnEstimationParticleTransmission1951}. It is published in the proceedings of a symposium held in  June 1949 but published in 1951 \cite{analysisu.s.MonteCarloMethod1951}. However, it is cited in 1950 in \cite{curtissSamplingMethodsApplied1950}.  It seems that the expression \textit{importance sampling} may have been used during the symposium. However, in this paper, there is no real explanation of the expression:

\begin{quote}
    For convenience, we refer to this second type of sampling where the neutrons are forced into the more important regions, as importance sampling.
\end{quote}

Nor is there a mention of Goertzel as the one who suggested the term as in \cite{kahnModificationMonteCarlo1950}.

Note that \cite{goertzelMonteCarloMethods1950} mentions an article with the same title but with the reference \textit{RAND Report R-148}. If the numbering of RAND reports is chronological, it would make it precede \cite{kahnStochasticMonteCarlo1949}, which has the reference \textit{R-163}. However we were not able to find any copy of this report online, and there is no guarantee that this draft contains the expression \textit{importance sampling} that may have been introduced afterwards in the final version.

\section{Conclusion}
It thus seems to us that the first time that \textit{importance sampling} was used was in \cite{kahnEstimationParticleTransmission1951} in 1949, maybe during the aforementioned June symposium. However, there is no mention of Goertzel found in this paper, while \cite{kahnModificationMonteCarlo1950} written afterwards, mentions Goertzel as the one who suggested the name. Did Kahn mention it only in \cite{kahnModificationMonteCarlo1950} because he knew the paper would be published before \cite{kahnEstimationParticleTransmission1951}, even though the conference had been held later ?

\bibliographystyle{alpha}
\bibliography{is_biblio}

\end{document}